\documentclass[12pt]{article} 

\usepackage[utf8]{inputenc} 
\usepackage[margin=2.54cm]{geometry}
\usepackage[centertags]{amsmath}
\usepackage{amsfonts}
\usepackage{newlfont}
\usepackage{amsthm,amsmath,amssymb,amscd,verbatim,epsfig} 
\usepackage{makeidx}
\usepackage{layout}
\usepackage{color}

\theoremstyle{definition}
\newtheorem{thm}{Theorem}[section]

\theoremstyle{definition}

\theoremstyle{remark}

\numberwithin{equation}{section}

\begin{document}

\title{The Best Bounds for Range Type Statistics}

\author{Tsung-Lin Cheng\footnote{Correspondent author: tlcheng@cc.ncue.edu.tw;
 National Changhua University of Education, Taiwan}\quad and Chin-Yuan Hu\footnote{Retired Professor, National Changhua University of Education, Taiwan} }

\maketitle
\begin{abstract}
In this paper, we obtain the upper and lower bounds for two inequalities related to the range statistics. The first one is concerning the one-variable case and the second one is
about the bivariate case.\\
\par
{\bf\emph{Keyword: Range Statistics, Inequality}}
\end{abstract}
\newpage
\section{Introduction}
For a long-memory time series $\{X_t,t=0,1,2,\cdots\}$, to estimate the so-called ``Hurst parameter'' $H$, \cite{H} used the adjusted rescaled range statistic (ARRS) defined by
\begin{equation}\label{RS}
  \frac{\max\limits_{1\le j\le n}\sum\limits_{i=1}^j(X_j-\bar{X})-\min\limits_{1\le j\le n}\sum\limits_{i=1}^j(X_j-\bar{X})}{\Big\{\frac1{n}\sum\limits_{i=1}^n (X_i-\bar{X})^2\Big\}^{1/2}}
\end{equation}
The asymptotic statistical properties of the ARRS was studied early on by \cite{Mand}. Its limit distribution is nonstandard but like the long-tailed distribution.
Let $x_1\le x_2\le \cdots\le x_n$ be a real sequence with elements not all equal, $n\ge 2$. \cite{DHP} provided tables of certain upper and lower percentage points of ${\cal R}/{\cal S}$ for samples of n observations from a single normal population, where
 the range statistic is defined by ${\cal R}_n=x_n-x_1>0$ and the sample variance is defined by
${\cal S}_n^2=\frac1{n-1}\sum\limits_{i=1}^n (x_i-\bar{x}_n)^2=\frac1{n(n-1)}\sum\limits_{i<j}|x_i-x_j|^2$, where $\bar{x}_n=\frac1n \sum\limits_{i=1}^n x_i$ stands for
the sample mean. Later,
 \cite{Thom} obtained the optimal upper and lower bounds for ${\cal R_n^2/S_n^2}$, namely
\begin{equation}
2(n-1)\ge {\cal R_n^2/S_n^2}
\ge
\left\{\begin{array}{ll}
                 \frac{4(n-1)}{n}, & \mbox{if $n$ is even} \\  
                 &\\
                 \frac{4n}{n+1}, & \mbox{if $n$ is odd}   
                               \end{array} \right.
 \end{equation}
 which is equivalent to
 \begin{equation}\label{thom57}
\frac{n}{2}\le\frac{\sum\limits_{i<j}|x_i-x_j|^2}{\cal R_n^2}
\le
\left\{\begin{array}{ll}
                 \frac{n^2}{4},\quad & \mbox{if $n$ is even} \\  
                 &\\
                 \frac{n^2-1}{4},\quad & \mbox{if $n$ is odd}.   
                               \end{array} \right.
 \end{equation}
 It is clear that the following bounds hold for all $p>0$, $n=2,3,\cdots$
 \[1\le \frac{\sum\limits_{i<j}|x_i-x_j|^p}{\cal R_n^p}\le \frac{n(n-1)}2. \]
 \cite{Thom} obtained the best bound for $p=2$ while the proof was not available there.
In this paper, we prove the best bounds for the case as $p\ge 1$.
In particular, as $p=1$, we are dealing with the Gini Mean Difference which is defined by
\begin{eqnarray}
  G_n &=& \frac1{n(n-1)}\sum\limits_{i=1}^n\sum\limits_{j=1}^n |x_i-x_j|= \frac2{n(n-1)}\sum\limits_{1\le i<j\le n} |x_{(i)}-x_{(j)}|\\
   &=&\frac4{n(n-1)}\sum\limits_{i=1}^n(i-\frac{(n+1)}2)x_{(i)}.
\end{eqnarray}
(see \cite{DaNa}).
However, the case when $0<p<1$ is still unknown.
In the following discussion, we will adopt a step-by-step transformation method to achieve the best upper and lower bounds for the inequalities of the form (\ref{TH1}) and (\ref{TH2}) as shown in Theorem 1 and 2.
\section{Main Results}

\begin{thm}\label{TH1}{\it
For positive integer $n\ge 2$, let $x_1\le \cdots\le x_n$ be real numbers that are not all equal.
If $p\ge 1$ is a positive real numbers, then

  \begin{equation}
\frac{n-2}{2^{p-1}}+1\le \frac{\sum\limits_{i<j}|x_i-x_j|^p}{|x_n-x_1|^p}
\le
\left\{\begin{array}{ll}
                 \frac{n^2}{4},\quad & \mbox{if $n$ is even} \\  
                 &\\
                 \frac{n^2-1}{4},\quad & \mbox{if $n$ is odd}.   
                               \end{array} \right.
 \end{equation}

\begin{itemize}

\item[a.] The upper bound is achieved\\
           \begin{description}
              \item[i.] as $n=2m, m\ge 1$, $x_1=\cdots=x_m<x_{m+1}=\cdots=x_n$, and
              \item[ii.] as $n=2m+1$, $x_1=\cdots=x_{m+1}<x_{m+2}=\cdots=x_n$ \\
              or $x_1=\cdots=x_m<x_{m+1}=\cdots=x_n$.
           \end{description}
\item[b] The lower bound is achieved when $x_1<x_2=\cdots=x_{n-1}<x_n$ and\\ $x_2=\cdots=x_{n-1}=\frac{x_n+x_1}2$.
\end{itemize}
}
\end{thm}
\begin{proof}

The term $\frac{\sum\limits_{i<j}|x_i-x_j|^p}{|x_n-x_1|^p}$ remains unchanged after scaling. Therefore, without loss of generality, we can assume the range statistic
${\cal R}_n=x_n-x_1$ to be a constant, e.g. ${\cal R}_n=1$. \\
\vskip1cm
({\bf\it a-i}).
First, as $n=2m$, $m\ge 1$, we prove
 \[\frac{\sum\limits_{i<j}|x_i-x_j|^p}{|x_n-x_1|^p}\le \frac{n^2}4.\]
We will show that under the condition $x_1=\cdots=x_m<x_{m+1}=\cdots=x_n$ the best upper bound $\frac{n^2}4$ is achieved. In this case, the best solution is
\[\frac{\sum\limits_{i<j}|x_i-x_j|^p}{|x_n-x_1|^p}=mm=\frac{n^2}4.\]
Consider the case as $x_1=\cdots=x_{n_1}, n_1\ge 1$ and $x_{n-n_2+1}=\cdots=x_n$, where $n_1\ge 1, n_2\ge 1$ are two positive integers.
Note that it doesn't require that $n_1+n_2=n$.

Firstly, we consider the case when $n_1<n_2$.
\begin{itemize}
\item[1.] When $n_1+n_2=2m=n$, if we move one point from the upper side to the lower side, namely the lower side contains $n_1+1$ points and
the upper side contains $n_2-1$ points, and then consider the change of the quantity  $\sum\limits_{i<j}|x_i-x_j|^p$.
In fact, the change of the quantity $\sum\limits_{i<j}|x_i-x_j|^p$ is
\[(n_1+1)(n_2-1)(x_n-x_1)^p-n_1n_2(x_n-x_1)^p,\]
which is  nonnegative. We can also observe that, the sum attains its maximum as $n_1=n_2=m$.

\item[2.] When $n_1<n_2$ and $n_1+n_2<n$, set $a=x_{n_1+1}-x_{n_1}>0$ and $b=x_{n}-x_{n_1+1}>0$. It is clear that ${\cal R}_n=x_n-x_1=a+b>0$.
If we add $x_{n_1+1}$ to the lower group, i.e. $x_1=x_2=\cdots=x_{n_1}=x_{n_1+1}$, then the change of the quantity  $\sum\limits_{i<j}|x_i-x_j|^p$ becomes
\begin{eqnarray*}
   n_2[(a+b)^p-b^p]-n_1a^p+\sum\limits_{j=n_1+2}^{n_2-1} (x_j-x_1)^p&> &  n_2[(a+b)^p-b^p]-n_1a^p  \\
     &\ge& n_1[(a+b)^p-b^p-a^p] >0
\end{eqnarray*}
Continue in this fashion until both sides having equal number of elements,i.e. $n_1=n_2$. Note that the number $n_1$ of the size of the lower side is not the older one but equals $n_2$. Since $n=2m$ is even, if there are no distinct values between $x_{n_1}$ and $x_{n-n_2+1}$, then $n_1=n_2=m$.

If there are some other distinct values between $x_{n_1}$ and $x_{n-n_2+1}$, then the number of those set must also be even.
We might assume that $a=x_{n_1+1}-x_1$(=$x_{n_2+1}-x_1)$, $b=x_{n-n_2}-x_{n_1+1}$ and $c=x_n-x_{n-n_2}$. It is clear that
$x_n-x_1=a+b+c$.

In this case, we may add $x_{n_1+1}$ to the lower side and $x_{n-n_2}$ to the upper side.
The the quantity  $\sum\limits_{i<j}|x_i-x_j|^p$ increases at least
\begin{eqnarray*}
   && n_2\bigg\{[(a+b+c)^p-(b+c)^p]+[(a+b+c)^p-(a+b)^p] \\
   && -(a^p+c^p)\bigg\}\\
   &>&0.
\end{eqnarray*}
Continue in this way until both sides have size $n_1=n_2=m$. We note that the value of $\sum\limits_{i<j}|x_i-x_j|^p$ increases after every reallocation.
We may also apply a parallel algorithm to the case as $n_1>n_2$ and will also have the same conclusion that $n_1=n_2=m$ is the optimal case.
\vskip1cm
({\bf\it a-ii}). When $n=2m+1$ is odd, we will show that
\[\frac{\sum\limits_{i<j}|x_i-x_j|^p}{|x_n-x_1|^p}\le \frac{n^2-1}4,\]
 which is equivalent to show that the optimal allocation occurs when either $n_1=m+1, n_2=m$ or $n_1=m, n_2=m+1$.
 We assume the lower side has $n_1$ elements and the upper side has $n_2$ elements, $n_1+n_2\le n$, i.e.
 \begin{eqnarray*}
   x_1 &=& x_2=\cdots =x_{n_1},\quad\hbox{ and } \\
   x_{n-n_2+1}&=&x_{n-n_2+2}=\cdots=x_n.
 \end{eqnarray*}
 We may prove this case in a similar fashion as above by considering the following three cases: $n_1<n_2$, $n_1=n_2$, and $n_1>n_2$.
 The first and third cases can be dealt with by removing one point to the smaller side until both sides have equal number of elements. In this case,
 one will achieve the condition when only one element is left between the lower sides and the upper side. Either adding the left point to the lower side or to the upper
 side is the optimal allocation.

 In the second case, we add one point to either side in each step until there is only one point is left. Similarly, either adding the left point to the lower side or to the upper
 side is the optimal allocation.
 \vskip1cm
 ({\bf\it b}). Now we prove the lower bound of the inequality. We just have to show that the optimal solution occurs when $x_1<x_2=\cdots=x_{n-1}<x_n$ and
 $x_2=\cdots=x_{n-1}=\frac{x_n+x_1}2$. Assume the lower side has $n_1$ points and the upper side has $n_2$ points, namely
  $x_1=\cdots=x_{n_1}, n_1\ge 1$ and $x_{n-n_2+1}=\cdots=x_n$, where $n_1\ge 1, n_2\ge 1$ are two positive integers.

 Consider the following three cases: $n_1<n_2$, $n_1=n_2$, and $n_1>n_2$. In the case as $n_1<n_2$, we adopt the following steps:\\
 1. As $n_1=1$, it requires no adjustment.\\
 2. As $n_1>1$, we move a point from the lower side to the middle part. Denote the distance between the lower side and the new location of the moved point by $a$ and
 that between the upper side and the new location of the moved point by $b$. It is clear that $x_n-x_1=a+b$.
 After such a reallocation, the quantity $\sum\limits_{i<j}|x_i-x_j|^p$ increases by the amount $(n_1-1)a^p$ and decreases at least by the amount $n_2[(a+b)^p-b^p]$.

 Note that
  \begin{eqnarray*}
   &&n_2[(a+b)^p-b^p]-(n_1-1)a^p   \\
   &\ge& (n_1-1)[(a+b)^p-b^p-a^p]>0.
 \end{eqnarray*}
 Therefore, $\sum\limits_{i<j}|x_i-x_j|^p$ decreases after such a reallocation.

 Repeat the above steps until $n_1=1$. In this case, $x_1<x_2\le \cdots<x_{n-n_2+1}=\cdots=x_n$, $n_2\ge 1$.
 If we move a point from the upper side to the middle part where it is closest to the leftist point of the upper side, we may also prove that
 $\sum\limits_{i<j}|x_i-x_j|^p$ decreases after such a reallocation. Repeat in this fashion until the upper side contains a single point, i.e.
 \[x_1<x_2\le \cdots\le x_{n-1}<x_n.\]
 $\sum\limits_{i<j}|x_i-x_j|^p$ will achieves its lower bound.

 In the case as $n_1=n_2$, we can also adopt the above procedures until both of the lower side and the upper side contains only a single point.
 If $n_1=n_2=1$, then stop reallocation. If $n_1=n_2>1$, then we move one point from both sides to the nearest location. Such a reallocation is
 the reversed procedures for proving the upper bound of the inequality. Repeat these steps until each side contains only a single point, i.e.
 \[x_1<x_2\le \cdots\le x_{n-1}<x_n.\]
 We can also deal with the case as $n_1>n_2$ in a similar fashion to the proof of the upper bound of the inequality.
  Next, we will show that the best allocation is when $x_1<x_2=\cdots=x_{n-1}<x_n$ and $x_2=\cdots=x_{n-1}=\frac{x_n+x_1}2$.

Suppose that $O$ is the mid-point $\frac{x_n+x_1}2$ of the range. The manner of removing points are to allocate these $n-2$ points in the
middle part as close as possible, to $O$. There are three possible cases: (i) the points of the middle part lie on the right side of $O$, (ii)  the points of the middle part lie on the right side of $O$, (iii) the points of the middle part lie on both sides of $O$.
In the case (i), suppose that among the $n-2$ points in the middle part, there are $m_1$ points $x_{n-m_1}=\cdots=x_{n-1}$ which are the nearest points from $x_n$. Assume in addition
that the distance between $x_{n-1}$ and $x_{n}$ is $a=x_n-x_{n-1}\ge 0$, the distance between $x_{n-m_1}$ and $x_{n-m_1-1}$ is $b=x_{n-m_1}-x_{n-m_1-1}\ge 0$, and the distance between
$x_{n-m_1-1}$ and $x_1$ is $c=x_{n-m_1-1}-x_1\ge a$. We then have $x_n-x_1=a+b+c$. If we move these $m_1$ points to the location of $x_{n-m_1-1}$, the quantity $\sum\limits_{i<j}|x_i-x_j|^p$ increases by $m_1[(a+b)^p-a^p]$ and decreases at least by $m_1[(c+b)^p-c^p]$. Note that the function $f(x)=(x+b)^p-x^p$ is an increasing function of $x\ge 0$ and $c\ge a$. Therefore,
\[m_1[(c+b)^p-c^p]-m_1[(a+b)^p-a^p]>0,\] and $\sum\limits_{i<j}|x_i-x_j|^p$ decreases at least by a positive amount after such a reallocation.

Repeating the above steps until all $n-2$ points of the middle part being moved to the center $\frac{x_1+x_n}2$, we achieve the allocation
\[x_1<x_2=x_3=\cdots=x_{n-1}=\frac{x_1+x_n}2<x_n.\]

Next, we consider case (ii). When all $n-2$ points of the middle part are located on the right of the center $O$, we can treat it as a dual case of case (i). Similarly, we
can prove that when all of the $n-2$ points are located at $O$, the quantity $\sum\limits_{i<j}|x_i-x_j|^p$ achieves its lower bound.

In case (iii), we may first move all points on the right side of $O$ to the center $O$ and then move all points on the left side of $O$ to $O$. We can also prove
that the lower bound is achieved when $x_1<x_2=\cdots=x_{n-1}=\frac{x_n+x_1}2<x_n$. Since the allocation requires only finite steps, such an allocation is an optimal solution.
\end{itemize}
\end{proof}
\begin{thm}\label{TH2}{\it For $n\ge 2$, $x_1\le\cdots\le x_n$ are not all equal and $y_1\le\cdots\le y_n$ are not all equal, we have
\begin{equation}\label{cov}
\frac{1}{n}\le \frac{\sum\limits_{i=1}^n(x_i-\bar{x})(y_i-\bar{y})}{(x_n-x_1)(y_n-y_1)}
\le
\left\{\begin{array}{ll}
                 \frac{n}{4},\quad & \mbox{if $n$ is even} \\  
                 &\\
                 \frac{n^2-1}{4n},\quad & \mbox{if $n$ is odd}.   
                               \end{array} \right.
 \end{equation}
\begin{itemize}
  \item[a.] Upper bound:
  \begin{description}
     \item[i.] In the case as $n=2m$, $m\ge 1$, the upper bound is achieved when
               \begin{eqnarray*}
                 x_1 &=& \cdots=x_m<x_{m+1}=\cdots=x_n\quad\hbox{ and } \\
                 y_1 &=& \cdots=y_m<y_{m+1}=\cdots=y_n.
               \end{eqnarray*}
     \item[ii.] In the case as $n=2m+1$, $m\ge 1$, the upper bound is achieved when
               \begin{eqnarray*}
                 x_1 &=& \cdots=x_m<x_{m+1}=\cdots=x_n\quad\hbox{ and } \\
                 y_1 &=& \cdots=y_m<y_{m+1}=\cdots=y_n.
               \end{eqnarray*}
               or
               \begin{eqnarray*}
                 x_1 &=& \cdots=x_{m+1}<x_{m+2}=\cdots=x_n\quad\hbox{ and } \\
                 y_1 &=& \cdots=y_{m+1}<y_{m+2}=\cdots=y_n.
               \end{eqnarray*}
     \end{description}
  \item b. The lower bound is achieved as
           \begin{eqnarray*}
             x_1<x_2=\cdots=x_n &\hbox{ and }& y_1=y_2=\cdots=y_{n-1}<y_n, \quad\hbox{ or } \\
             x_1=x_2=\cdots=x_{n-1}<x_n &\hbox{ and }& y_1<y_2=\cdots=y_n.
           \end{eqnarray*}
\end{itemize}
}
\end{thm}
\begin{proof}
Without loss of generality, we may assume that $x_1=y_1=0$ and $x_n=y_n=1$ since scaling and shifting don't change the value of (\ref{cov}). The quantity of (\ref{cov}) then reduces to
$s_{xy}\equiv \sum\limits_{i=1}^n(x_i-\bar{x})(y_i-\bar{y})=\sum\limits_{i=1}^nx_iy_i-n\bar{x}\bar{y}$. The values of $x's$ and $y's$ are not all equal. The number of all possible allocations is $4^{n-1}$, which is an NP-hard problem for computer science.
\begin{itemize}
    \item[a.] \underline{\bf Proof of the Upper Bound }
    Since $0=y_1\le y_2\le\cdots\le y_n=1$ are not all equal, $0<\bar{y}<1$. There exists a positive $k$ with $2\le k\le n$ such that $y_{k-1}\le \bar{y}\le y_k$. Now for such $k$, we replace $x_k$ by $x'_k=x_k+a$ for some undetermined constant $a$, and other $x_i, i\ne k$ remain unchanged. we denote the new sequence by $\{x'_i,i=1,2,\cdots,n\}$. In this case, the new mean becomes $\bar{x'}=\bar{x}+\frac{a}{n}$ and the corresponding $s_{x'y}$ becomes
    \begin{eqnarray*}
    s_{x'y}&\equiv & \sum\limits_{i=1}^n (x'_i-\bar{x'})(y_i-\bar{y})=\sum\limits_{i=1}^n x'_iy_i-n\bar{x'}\bar{y}\\
    &=& \sum\limits_{i\ne k}x_iy_i+(x_k+a)y_k-n(\bar{x}+\frac{a}{n})\bar{y}\\
    &=& s_{xy}+a(y_k-\bar{y}).
    \end{eqnarray*}
    Due to the fact that $y_k\ge \bar{y}$, if $a\ge 0$, then $s_{x'y}\ge s_{xy}$.

    Now, we move $x_k$ to the same position as $x_{k+1}$, i.e. $a=x_{k+1}-x_k$, then the quantity $s_{xy}$ increases.
    Since $y_{k+1}\ge y_k\ge\bar{y}$, we may define a new $x's$ sequence by setting only $x'_{k+1}=x_{k+1}+a$ as above and other $x_i,i\ne k+1$, remain unchanged.
    We may repeat this procedure until
    $x_k=x_{k+1}=\cdots=x_n=1.$ In this case $x_k\ge \bar{x}$. In this case, we may define a new $y's$ sequence by setting $y'_k=y_k+b$ for some undetermined constant $b$.
    For the new $x's$ and $y's$ sequences, observe that
     \begin{eqnarray*}
    s_{xy'}&\equiv & \sum\limits_{i=1}^n (x_i-\bar{x})(y'_i-\bar{y'})\\
        &=& s_{xy}+b(x_k-\bar{x}).
    \end{eqnarray*}
    If $b\ge 0$, then the quantity $s_{xy'}$ for the new sequences increases. Similar to the above argument, we may choose $b=y_{k+1}-y_k$. Repeat this procedure until we have
    $y_k=y_{k+1}=\cdots=y_n=1$.

    Note that, $y_{k-1}\le \bar{y}\le y_k$, we can adopt a reversed procedure to allocate the lower side by letting $x'_{k-1}=x_{k-1}-a$ for some undetermined constant $a$ and other terms be unchanged. In this case, the quantity $s'_{xy}$ becomes $s'_{xy}=s_{xy}-a(y_{k-1}-\bar{y})$. Similarly, $s'_{xy}\ge s_{xy}$ if $a\ge 0$. Now we may let $a=x_{k-1}-x_{k-2}$. In this case, $y_{k-2}\le y_{k-1}\le \bar{y}$, we may repeat the above procedure until we have $x_1=\cdots=x_{k-1}=0$. In a similar way to the above procedures for allocating $y's$, we may finally achieve the condition as $y_1=\cdots=y_{k-1}=0$. Therefore, the upper bound can be achieved as, for $\ell=k-1$, $1\le\ell\le n-1$,
    \begin{eqnarray*}
                 0=x_1 &=& \cdots=x_{\ell}<x_{\ell+1}=\cdots=x_n=1\quad\hbox{ and } \\
                 0=y_1 &=& \cdots=y_{\ell}<y_{\ell+1}=\cdots=y_n=1,\quad\hbox{ for some }1\le\ell\le n-1.
               \end{eqnarray*}
    In this case, $\sum\limits_{i=1}^n x_i=n-\ell=\sum\limits_{i=1}^n y_i$, $\sum\limits_{i=1}^n x_iy_i=n-\ell$, and
    \begin{eqnarray*}
    s_{xy}&=& \sum\limits_{i=1}^n (x_i-\bar{x})(y_i-\bar{y})=\sum\limits_{i=1}^n x_iy_i-n\bar{x}\bar{y}\\
    &=&(n-\ell)-\frac{(n-\ell)^2}{n}\\
    &\le& \left\{\begin{array}{ll}
                 \frac{n}{4},\quad & \mbox{if $n$ is even} \\  
                 &\\
                 \frac{n^2-1}{4n},\quad & \mbox{if $n$ is odd}.   
                               \end{array} \right. .
    \end{eqnarray*}

    \begin{description}
    \item[(i)] When $n=2m$, $m\ge 1$, the optimal bound occurs as $\ell=m$.
   \item[(ii)]  When $n=2m+1$, $m\ge 1$, the optimal bound occurs as $\ell=m$ or $m+1$.

    \end{description}
    \item[b.] \underline{\bf Lower Bound}
Due to the condition $0\le y_1\le\cdots\le y_n=1$, we have $0<\bar{y}<1$. There exists $2\le k\le n$ such that $y_{k-1}\le \bar{y}\le y_k$.

\begin{description}
\item[(1)] As $\bar{y}\le y_k$, For such $k$, we set $x'_k=x_k-a$ for some undetermined $a$ and make the other $x_i's$, $i\ne k$, remain unchanged.  The new average $\bar{x'}=\bar{x}-\frac{a}{n}$ and the cross term $s_{x'y}=s_{xy}-a(y_k-\bar{y})$. If $a\ge 0$, then $s_{x'y}\le s_{xy}$. We may set $a=x_k-x_{k-1}$.
It is clear that $y_{k+1}\ge \bar{y}$. We may set $x'_{k+1}=x_{k+1}-a$  for some undetermined $a$ and continue in the above way until
$x_{k-1}=x_k=x_{k+1}=\cdots=x_{n-1}\le x_n=1$.
\item[(2)] On the other hand, as $y_{k-1}\le \bar{y}$, set $x'_{k-1}=x_{k-1}+a$ for some undetermined $a$ and make the other $x_i$-terms remain unchanged.
 Then $\bar{x'}=\bar{x}+\frac{a}{n}$ and the cross term $s_{x'y}=s_{xy}+a(y_{k-1}-\bar{y})$. If $a\ge 0$, then $s_{x'y}\le s_{xy}$. We may choose $a=x_k-x_{k-1}$.
  Note that $y_{k-2}\le \bar{y}$. In a similar fashion, we may reallocating $x_{k-2}, x_{k-3},\cdots$, until $x_1\le x_2=\cdots=x_{k-1}$.
\end{description}
Combining (1) and (2), we obtain that as $0=x_1\le x_2=\cdots=x_{n-1}\le x_n=1$. Similarly, we can reallocating the $y_i's$ until $0=y_1\le y_2\cdots=y_{n-1}\le y_n=1$. Next, we will
show that the lower bound is achieved as
           \begin{eqnarray*}
             x_1<x_2=\cdots=x_n &\hbox{ and }& y_1=y_2=\cdots=y_{n-1}<y_n, \quad\hbox{ or } \\
             x_1=x_2=\cdots=x_{n-1}<x_n &\hbox{ and }& y_1<y_2=\cdots=y_n.
           \end{eqnarray*}
Set $x_2=\cdots=x_{n-1}=a$ and $y_2=\cdots=y_{n-1}=b$, where $0\le a,b\le 1$.
In this case, $\sum\limits_{i=1}^n x_i=a(n-2)+1$ and $\sum\limits_{i=1}^ny_i=b(n-2)+1$, and $\sum\limits_{i=1}^n x_iy_i=ab(n-2)+1$.
As a result
   \begin{eqnarray*}
    s_{xy}&=&\sum\limits_{i=1}^n x_iy_i-n\bar{x}\bar{y}\\
    &=&\Big[ab(n-2)+1\Big]-n\Big(\frac{a(n-2)+1}{n}\Big)\Big(\frac{b(n-2)+1}{n}\Big)\\
    &=& \Big(\frac{n-2}{n}\Big)[ab+(1-a)(1-b)]+\frac1{n}\\
    &\ge& \frac1{n},
    \end{eqnarray*}
    and the equality holds if $a=1,b=0$, or, $a=0,b=1$.
  \end{itemize}
\end{proof}
\vskip1cm
{\bf Remark.} We may check the optimal cases as below:
  \begin{itemize}
    \item[a.] \underline{\bf Upper Bound }
    \begin{description}
    \item[(i)] When $n=2m$, $m\ge 1$,
    \begin{eqnarray*}
                 0=x_1 &=& \cdots=x_m<x_{m+1}=\cdots=x_n=1\quad\hbox{ and } \\
                 0=y_1 &=& \cdots=y_m<y_{m+1}=\cdots=y_n=1,
               \end{eqnarray*}
               as a result $\sum\limits_{i=1}^nx_iy_i-n\bar{x}\bar{y}=m-\frac{m^2}n=\frac{n}4$.
     \item[(ii)]  When $n=2m+1$, $m\ge 1$,
       \begin{eqnarray*}
                 0=x_1 &=& \cdots=x_m<x_{m+1}=\cdots=x_n=1\quad\hbox{ and } \\
                 0=y_1 &=& \cdots=y_m<y_{m+1}=\cdots=y_n=1,
               \end{eqnarray*}
           we have $\sum\limits_{i=1}^nx_iy_i-n\bar{x}\bar{y}=m+1-n(\frac{m+1}{n})^2=\frac{n^2-1}{4n}$.
    \end{description}
    \item[b.] \underline{\bf Lower Bound} Let
   $0=x_1<x_2=\cdots==x_n=1$ and $0=y_1=\cdots=y_{n-1}<y_n=1$. we have
   \[\sum\limits_{i=1}^nx_iy_i-n\bar{x}\bar{y}=\frac1{n}.\]

  \end{itemize}


\end{document}